\documentclass[a4paper,12pt,twoside]{article}
\usepackage[T1]{fontenc}
\usepackage{lmodern}
\usepackage{fullpage}
\usepackage{amssymb,amsmath,amsthm,mathtools}
\usepackage{enumitem}
\usepackage{needspace}
\usepackage[unicode,hidelinks]{hyperref}
\usepackage{url}

\newtheorem{theorem}{Theorem}

\newtheorem{lemma}{Lemma}

\newcommand{\N}{\mathbb N}
\newcommand{\Z}{\mathbb Z}
\newcommand{\ld}{\underline d}
\newcommand{\ord}{\operatorname{ord}}
\newcommand{\ordR}{\operatorname{ord}_{R}}

\hypersetup{pdftitle={On a Problem of Erd\H{o}s and Graham Related to Restricted Order},pdfauthor={Shi-Qiang Chen, Quan-Hui Yang, Lilu Zhao}}
\begin{document}
\title{\Large On a problem of Erd\H{o}s and Graham
related to restricted order\footnote{The first author is supported by the National Natural Science Foundation of China (Grant No. 12301003), the Anhui Provincial Natural Science Foundation (Grant No. 2308085QA02) and the University Natural Science Research Project of Anhui Province (Grant No. 2022AH050171). The second author is supported by the National Natural Science Foundation of China (Grant No. 12371005 and 12671008). The third author is supported by the National Key Research and
Development Program of China (Grant No.~2021YFA1000701) and the
National Natural Science Foundation of China (Grant No.~12471088). }}

\author{\large Shi-Qiang Chen$^{a}$,
\large ~~~Quan-Hui Yang$^{b}$,~~~Lilu Zhao$^{c}$,\thanks{E-mail addresses: csq20180327@163.com (S.-Q. Chen), yangquanhui01@163.com (Q.-H. Yang), zhaolilu@ustc.edu.cn (L. Zhao).}}
\date{} \maketitle
\par\vspace{-1.5em}
\begingroup
\centering
\small
\setlength{\parskip}{0pt}
\setlength{\baselineskip}{13pt}
$a$. School of Mathematics and Statistics,\\
Anhui Normal University, Wuhu 241002, China\par
\vspace{3pt}
$b$. School of Mathematical Sciences and Ministry of Education Key Laboratory for NSLSCS,\\
Nanjing Normal University, Nanjing 210023, China\par
\vspace{3pt}
$c$. School of Mathematical Sciences,\\
University of Science and Technology of China, Hefei 230026, China\par
\endgroup
\vspace{1em}

\par\noindent {\bf Abstract.}\quad
Let $A$ be a set of nonnegative integers. We prove that if
$A\setminus F$ has a finite order for every
finite subset $F\subseteq A$, then $A$ has finite restricted
order. This resolves a problem posed by Erd\H{o}s and Graham in
1980. We also establish a local to global principle to characterize restricted asymptotic bases of finite order.

\medskip
\noindent{\bf Keywords:} asymptotic basis; order; restricted order; local to global principle.

\medskip
\noindent 2020 {\it Mathematics Subject Classification}: 11B13; 11B05.

\section{Introduction}
Let $\mathbb N$ denote the set of all nonnegative integers. For a set $X$, let $|X|$ denote its
cardinality, with $|X|=\infty$ when $X$ is infinite.
For $X\subseteq\mathbb N$ and an integer $q\geq1$, write
\[
\begin{aligned}
X\bmod q&=\{x+q\mathbb Z:x\in X\}\subseteq\mathbb Z/q\mathbb Z.
\end{aligned}
\]
For $X\subseteq\mathbb N$, we define the lower density
\[
\begin{aligned}
\ld(X)&=\liminf_{x\to\infty}\frac{|X\cap[1,x]|}{x}.
\end{aligned}
\]
For $A\subseteq\mathbb N$ and an integer $j\geq0$, define $S_j(A)=\bigcup\limits_{i=0}^j iA$,
$$R_j(A)=\left\{\sum_{a\in F}a:F\subseteq A,\ |F|\leq j\right\},$$
and $\Sigma(A)=\bigcup_{j\geq0}R_j(A)$. Here $iA$ is the set of sums of $i$ elements of $A$, where repetitions
are allowed, and $0A=\{0\}$. The summands in $R_j(A)$ are distinct, and the empty sum
is zero. For $X,Y\subseteq\mathbb N$, we write $X\sim Y$ if both
$X\setminus Y$ and $Y\setminus X$ are finite. Following Erd\H{o}s and Graham \cite{EG1980}, a set $A$ is said to be an
 asymptotic basis of finite order (or a basis in short) if $S_h(A)\sim\mathbb N$ for some integer $h\geq1$,
and a  restricted asymptotic basis of finite order if $R_h(A)\sim\mathbb N$ for some
integer $h\geq1$.

The ordinary order and the restricted order of $A$ are defined by
$$\ord(A)=\min\{h\geq1:S_h(A)\sim\mathbb N\},$$
$$\ordR(A)=\min\{h\geq1:R_h(A)\sim\mathbb N\}.$$
If the corresponding minimum does not exist, we set the order equal
to $\infty$.

The study of restricted order goes
back at least to Kelly~\cite{Kelly1957}, who proved that $\ordR(A)\leq4$ if
$\ord(A)\leq2$, and conjectured that $\ordR(A)\leq3$ if
$\ord(A)\leq2$. In fact, Kelly also proved that $\ordR(A)\leq3$ if $\ord(A)=2$
and $\ld(A)>0$. Kelly's conjecture was disproved by Hennecart~\cite{Hennecart2005}, who constructed $A$ satisfying $\ord(A)=2$ and $\ordR(A)=4$.
Hegyv\'ari, Hennecart and Plagne~\cite{HHP2003} studied the lower density of the finite restricted addition of a basis.

In general, an
asymptotic basis of finite order need not have finite
restricted order. For example, if $h\geq3$ and
$A=\{1\}\cup h\mathbb N$, then $\ord(A)\leq h$, but it does not have finite
restricted order. Moreover, $A\setminus\{1\}$ is not an asymptotic basis of finite order. This is observed by Bateman (see~\cite{Bloom338,EG1980}). This leads to the following question
whether finite restricted order follows if the finite ordinary order property
remains after every finite deletion.

In 1980, Erd\H{o}s and Graham~\cite{EG1980} posed the following problem.
\medskip

\noindent{\bf Erd\H{o}s--Graham's Problem.}
Let $A\subseteq\mathbb N$. Is it true that if
$\ord(A\setminus F)<\infty$ for every finite subset $F\subseteq A$, then
 $\ordR(A)<\infty$?
\medskip

In this paper, we  resolve the above Erd\H{o}s--Graham problem affirmatively.
\begin{theorem}\label{thm1}
Let $A\subseteq\mathbb N$. Suppose that $\ord(A\setminus F)<\infty$ for every finite subset
$F\subseteq A$. Then $\ordR(A)<\infty$.
\end{theorem}

In order to prove Theorem \ref{thm1}, we developed a local to global principle to characterize whether a  basis $A$ is a restricted asymptotic basis of finite order.

For $A\subseteq \mathbb N$ and $q\geq 1$, we say $A$ is complete modulo $q$ if
\[
\Sigma(A)\bmod q=\mathbb Z/q\mathbb Z.
\]
If $A$ has a finite restricted order, then one immediately deduces that $A$ is complete modulo $q$ for every integer $q\geq 1$. We prove that this is indeed a sufficient condition as well provided that $A$ has a finite ordinary order. Precisely, we have the following.

\begin{theorem}[The local to global principle]\label{thm2}
Let $A\subseteq\mathbb N$ with $\ord(A)<\infty$.
Then $\ordR(A)<\infty$ if and only if $A$ is complete modulo $q$ for every integer $q\geq 1$.
\end{theorem}

Erd\H{o}s and Graham \cite[p.206]{EG1980} also asked the following:

\noindent {\it ``What are necessary and sufficient conditions on a basis $A$ to have restricted order?''}

\noindent Theorem \ref{thm2} provides an answer to this question.

In view of Theorem \ref{thm2}, it is important to study the local property of an asymptotic basis. For convenience, we say that $A\subseteq \mathbb N$ is {\it locally complete} if $A$ is complete modulo $q$ for every integer $q\geq 1$. Note that if $A$ is complete modulo $q$, then there exists $H:=H_q$ such that $$R_H(A)\bmod q=\mathbb Z/q\mathbb Z.$$
For $A\subseteq\mathbb N$, we say that $A$ is
\emph{uniformly locally complete} if there exists a positive integer $H$ (independent of $q$)
such that, for every integer $q\ge1$,
$R_H(A)\bmod q=\mathbb Z/q\mathbb Z$.

The next result shows that a locally complete basis is uniformly locally complete.
\begin{theorem}\label{thm3}
Let $A\subseteq\mathbb N$ with $\ord(A)<\infty$. Suppose that $A$ is locally complete. Then $A$ is uniformly locally complete.
\end{theorem}
This is reminiscent of the well-known fact that a continuous function on a closed interval is uniformly continuous.

For results on additive bases, restricted addition and subset sums,
see~\mbox{\cite{DG2001,GGHHLP2002,Plagne2008,SV2006}}. Recently, Veljjanoski  \cite{Veljjanoski2026} made partial progress on the Erd\H{o}s--Graham problem for sets with positive lower density.

In Section 2, we prepare some auxiliary results.  In Section 3, we study the local complete property of asymptotic bases and complete the proof of Theorem \ref{thm3}. Then we establish the local to global principle in Section 4. Finally, we prove Theorem \ref{thm1} in Section 5.

\section{Preparations}

If there are an integer $g\ge1$ and a set $T\subseteq\Z/g\Z$ such that
\[
 X\sim\{n\in\mathbb N:n\bmod g\in T\},
\]
then we say that $X$ is eventually periodic with period $g$.
\begin{lemma}[\cite{Kneser}]\label{lem1}
Let $t\ge2$, and let $X_1,\ldots,X_t\subseteq\mathbb N$ be infinite sets; the sets $X_i$ are not required to be distinct.
If
\[
 \ld(X_1+\cdots+X_t)
 <\liminf_{N\to\infty}\frac{\sum_{i=1}^t|X_i\cap[0,N]|}{N},
\]
then $X_1+\cdots+X_t$ is eventually periodic.
\end{lemma}

\begin{lemma}[\cite{HHP2003}]\label{lem2}
Let $A\subseteq\mathbb N $ and let $h\ge1$. If
$S_h(A)\sim\mathbb N $, then
\[
 \ld(R_h(A))>0.
\]
\end{lemma}
\begin{proof}
Let $r=\ord(A)$. Then $r\le h$. Write
\[
 r\times A
 =\{a_1+\cdots+a_r:a_i\in A,\ a_i\ne a_j\ (i\ne j)\}.
\]
Note that $r\times A\subseteq R_h(A)$. Hegyv\'ari, Hennecart and Plagne~\cite[Theorem 4]{HHP2003} proved that
\[
  \ld(r\times A)>0,
\]
and the desired conclusion of Lemma~\ref{lem2} follows.
\end{proof}

For $r\ge1$ and $n\ge1$, let $\rho_{2r}(n;A)$ be the largest
integer $q\ge0$ for which there are pairwise disjoint nonempty sets
$F_1,\ldots,F_q\subseteq A$ satisfying $|F_i|\le2r$, and
$\sum_{a\in F_i}a=n$ for $i=1,2,\ldots,q$.
 Here $q=0$ is allowed, and the maximum exists since
$A\cap[0,n]$ is finite. For $r,Q\ge1$, put
\[
 D_{r,Q}=\{n\ge1:\rho_{2r}(n;A)\ge Q\}.
\]
\begin{lemma}\label{lem3}
Let $r,Q,N\ge1$, and put $\mathcal X=R_r(A)\cap[1,N]$. Write
$M=|\mathcal X|$ and $$e=|R_{r-1}(A)\cap[0,N]|.$$
For each $x\in\mathcal X$, fix a set $F_x\subseteq A$ with $1\le|F_x|\le r$
such that $ \sum\limits_{a\in F_x}a=x$. Put
$$\mathcal P=\{(x,y)\in\mathcal X^2:F_x\cap F_y=\emptyset\},$$
and put
$$p_n=\left|\{(x,y)\in\mathcal P:x+y=n\}\right|~~(n\ge0).$$
Then
$$\sum_{n=0}^{2N}p_n\ge M^2-rMe,$$
and $p_n\le M$. Moreover, $\{n:p_n>4r(Q-1)e\}\subseteq D_{r,Q}$.
\end{lemma}
\begin{proof}
Fix $a\in A$. If $y\in\mathcal X$ and $a\in F_y$, then
\[
 y-a=\sum_{u\in F_y\setminus\{a\}}u
 \in R_{r-1}(A)\cap[0,N].
\]
Hence
\begin{equation}\label{eq4}
 \left|\{y\in\mathcal X:a\in F_y\}\right|
 \le |R_{r-1}(A)\cap[0,N]|=e.
\end{equation}

If $(x,y)\notin\mathcal P$, then $F_x\cap F_y\ne\emptyset$.
By \eqref{eq4},
\begin{eqnarray*}M^2-|\mathcal P|
&=&\left|\{(x,y)\in\mathcal X^2:F_x\cap F_y\ne\emptyset\}\right| \\
&\le &\sum_{x\in\mathcal X}\sum_{a\in F_x}
       \left|\{y\in\mathcal X:a\in F_y\}\right| \le e\sum_{x\in\mathcal X}|F_x| \le rMe.
\end{eqnarray*}
Since $x,y\in[1,N]$, we have $2\le x+y\le2N$, and hence
\[
 \sum_{n=0}^{2N}p_n=\sum_{n=2}^{2N}p_n
 =|\mathcal P|\ge M^2-rMe.
\]
Next, fix $n$. Then
\[
 p_n\le\left|\{x\in\mathcal X:n-x\in\mathcal X\}\right|\le M.
\]

Fix $n\notin D_{r,Q}$. Among the ordered pairs satisfying
$(x,y)\in\mathcal P$ and $x+y=n$,
choose a maximal family of ordered pairs $\mathcal M$ such that the corresponding sets
$F_x\cup F_y$ are pairwise disjoint. Put
\[
 V=\bigcup_{(x,y)\in\mathcal M}(F_x\cup F_y).
\]
Since $F_x\cap F_y=\emptyset$, we have
$\sum_{a\in F_x\cup F_y}a=x+y=n$ and $|F_x\cup F_y|\le2r$.
Since $n\notin D_{r,Q}$, it follows that $|\mathcal M|\le Q-1$, and
$$|V|\le2r(Q-1).$$
Since $\mathcal M$ is maximal,
every ordered pair satisfying
$(x,y)\in\mathcal P$ and $x+y=n$ must satisfy
\[
 (F_x\cup F_y)\cap V\ne\emptyset.
\]
By \eqref{eq4} again,
for each $a\in V$, we have
\[\left|\{(x,y)\in\mathcal P:x+y=n,\ a\in F_x\}\right|\le e\]
and
\[ \left|\{(x,y)\in\mathcal P:x+y=n,\ a\in F_y\}\right|\le e.\]
Therefore,
\[
 \begin{aligned}
 p_n&\le\sum_{a\in V}\bigl(
 \left|\{(x,y)\in\mathcal P:x+y=n,\ a\in F_x\}\right|\\
 &~~~+\left|\{(x,y)\in\mathcal P:x+y=n,\ a\in F_y\}\right|\bigr)\\
 &\le2e|V|\le4r(Q-1)e.
 \end{aligned}
\]
Thus, if $p_n>4r(Q-1)e$, then $n\in D_{r,Q}$.

This completes the proof of Lemma~\ref{lem3}.
\end{proof}

\begin{lemma}\label{lem4}
Let $A\subseteq \mathbb N$, and suppose that $\ld(R_h(A))>0$ for some $h\ge1$.
For every integer $b\ge0$, there exist integers $K,g,N_0\ge1$ and $a$ with
$0\le a<g$ such that, for every finite set $F\subseteq A$ with $|F|\le b$ and
every integer $n\ge N_0$ satisfying $n\equiv a\pmod g$, we have
\[
 n\in R_K(A\setminus F).
\]
Here $K,g,N_0,a$ are independent of $F$.
\end{lemma}
\begin{proof}
Choose $0<\alpha_h<\ld(R_h(A))$, and define the following quantities successively for $r=h,h-1,\ldots,1$:
\begin{equation}\label{eq14}
L_r=\left\lceil 8\alpha_r^{-2}\right\rceil,\quad
Q_r=b+\sum_{j=r+1}^{h}2jL_j+2r(L_r-1)+1,\quad
\alpha_{r-1}=\frac{\alpha_r^2}{64rQ_r}.
\end{equation}

%

%
Put $D_r=D_{r,Q_r}$ and $C_r=D_r\cup\{0\}$, and set
\[
 K=\sum_{r=1}^h2rL_r.
\]
For sufficiently large $N$, choose the smallest index $r$ such that
\[
 |R_r(A)\cap[1,N]|\ge\alpha_rN.
\]
Such an index $r$ exists because $\ld(R_h(A))>\alpha_h$.
Write
\[
 M=|R_r(A)\cap[1,N]|,\quad
 e=|R_{r-1}(A)\cap[0,N]|.
\]
First, take $N$ sufficiently large so that
\[
 N\ge\max_{0\le j\le h}\frac1{\alpha_j}.
\]
If $r\ge2$, then by the minimality of $r$, we have $|R_{r-1}(A)\cap[1,N]|<\alpha_{r-1}N$, and
$$e=1+|R_{r-1}(A)\cap[1,N]|
 <1+\alpha_{r-1}N\le2\alpha_{r-1}N.$$

If $r=1$, then $R_0(A)=\{0\}$, and hence
\[
 e=1\le2\alpha_0N.
\]
Therefore, whether $r=1$ or $r\ge2$, we have $M\ge\alpha_rN$, and $e\le2\alpha_{r-1}N$.

Apply Lemma~\ref{lem3} with $Q=Q_r$, and let $p_n$ be as defined there.
On the one hand,
\[
 \sum_{n=1}^{2N}p_n\ge M^2-rMe.
\]
On the other hand, for $n\in D_r$ we use $p_n\le M$,
and for $n\notin D_r$ we use $p_n\le4r(Q_r-1)e$. Thus
\[
 \begin{aligned}
 \sum_{n=1}^{2N}p_n
 &\le M|D_r\cap[1,2N]|+4r(Q_r-1)e(2N)\\
 &\le N|D_r\cap[1,2N]|+8r(Q_r-1)eN.
 \end{aligned}
\]
Combining these two estimates and substituting $M\ge\alpha_rN$,
$e\le2\alpha_{r-1}N$, and \eqref{eq14}, we obtain
\[
 \begin{aligned}
 |D_r\cap[1,2N]|
 &\ge\frac{M^2-rMe}{N}-8r(Q_r-1)e\ge\frac{M^2}{N}-r(8Q_r-7)e\\
 &\ge\alpha_r^2N-2r(8Q_r-7)\alpha_{r-1}N =\alpha_r^2N\left(1-\frac{8Q_r-7}{32Q_r}\right)\ge\frac34\alpha_r^2N.
 \end{aligned}
\]
Since $L_r\ge8\alpha_r^{-2}$ and $|D_r\cap[1,2N]|\ge\tfrac34\alpha_r^2N$, it follows that
$$L_r|C_r\cap[0,2N]|
 \ge L_r|D_r\cap[1,2N]| \ge6N,$$
and so
$$\sum_{j=1}^hL_j|C_j\cap[0,2N]| \ge6N.$$
First, if $X=2N$ is a sufficiently large even integer, then
\[
 \sum_{j=1}^hL_j|C_j\cap[0,X]|\ge6N=3X.
\]
If $X=2N+1$ is a sufficiently large odd integer, the monotonicity of the counting functions gives
\[
 \sum_{j=1}^hL_j|C_j\cap[0,X]|
 \ge\sum_{j=1}^hL_j|C_j\cap[0,2N]|\ge6N.
\]
Dividing the lower bounds in both cases by $X$ and letting $X\to\infty$, we obtain
\[
 \liminf_{X\to\infty}\frac{\sum_{r=1}^hL_r|C_r\cap[0,X]|}{X}
 \ge3.
\]
Let $I=\{r:D_r\text{ is infinite}\}$ and put
\[
 S=\sum_{r\in I}L_rC_r.
\]
Since $C_r$ is finite whenever $r\notin I$, and there are only finitely many indices $r$,
removing these terms does not affect the lower bound for the limit inferior. Hence
\[
 \begin{aligned}
 \liminf_{X\to\infty}
 \frac{\sum_{r\in I}L_r|C_r\cap[0,X]|}{X}&\ge3>1\geq\ld(S),\\
 \end{aligned}
\]
In particular, $I\ne\emptyset$. Since $0<\alpha_r<1$ and
$L_r\ge8$ for every $r$, we have
\[
\sum_{r\in I}L_r\ge2.
\]
All the summand sets $C_r$ with $r\in I$ are infinite. Hence
Lemma~\ref{lem1}, applied to $L_r$ copies of each $C_r$, shows that
$S$ is eventually periodic. Moreover, $0\in C_r$ for every $r$, so
$C_r\subseteq S$ whenever $r\in I$, and $S$ is infinite.
Thus there are integers $g,N_0\ge1$ and $a$ with $0\le a<g$ such that
\begin{equation}\label{eq16}
 [N_0,\infty)\cap(a+g\Z)\subseteq S.
\end{equation}
Fix a finite set $F\subseteq A$ with $|F|\le b$, and take any $s\in S$. By the definition of $S$,
there are elements $d_{r,\ell}\in C_r$, where $C_r=D_{r,Q_r}\cup\{0\}$, such that
$$s=\sum_{r\in I}\sum_{\ell=1}^{L_r}d_{r,\ell}.$$
Arrange the elements of $I$ in decreasing order as follows:
$$I=\{r_1,r_2,\ldots,r_t\},~~
 h\ge r_1>r_2>\cdots>r_t\ge1.$$
Write
$$\mathcal J=\{(r_i,\ell):1\le i\le t,\ 1\le\ell\le L_{r_i}\}.$$
Define an order on $\mathcal J$ by
\[
 (r_i,\ell)\prec(r_j,m)
\quad\text{if and only if}\quad
i<j\ \text{or both}\  i=j\text{ and }\ell<m.
\]

Following this order, we prove by induction on $k$ that the set for the $k$th index pair can be chosen so that all sets $U_{r,\ell}\subseteq A$ chosen so far satisfy
\begin{equation}\label{eq17}
d_{r,\ell} =\sum_{u\in U_{r,\ell}}u,~~|U_{r,\ell}| \le2r,
~~U_{r,\ell}\cap F =\emptyset,
~~U_{r,\ell}\cap U_{j,m} =\emptyset~((j,m)\prec(r,\ell)).
\end{equation}
Let $k\ge1$, and suppose that the sets for the first $k-1$ index pairs have been chosen so that \eqref{eq17} holds. When $k=1$, no sets have yet been chosen.
Let $(r,\ell)$ be the $k$th index pair. By the definition of $\prec$,
the sets already chosen are all the sets with index $j>r$, together with
$U_{r,1},\ldots,U_{r,\ell-1}$. Put
\begin{equation}\label{eq18}
 W_{r,\ell}=F\cup
\bigg( \bigcup_{\substack{j\in I\\j>r}}\bigcup_{m=1}^{L_j}U_{j,m}\bigg)
 \cup\left(\bigcup_{m=1}^{\ell-1}U_{r,m}\right).
\end{equation}
By the induction hypothesis, each previously chosen set $U_{j,m}$ satisfies $|U_{j,m}|\le2j$.
Thus, using $|F|\le b$, $1\le\ell\le L_r$, and \eqref{eq14}, we have
\[
\begin{aligned}
 |W_{r,\ell}|
 &\le |F|+
 \sum_{\substack{j\in I\\j>r}}\sum_{m=1}^{L_j}|U_{j,m}|
 +\sum_{m=1}^{\ell-1}|U_{r,m}|\\
 &\le |F|+
 \sum_{\substack{j\in I\\j>r}}\sum_{m=1}^{L_j}2j
 +\sum_{m=1}^{\ell-1}2r\\
 &=|F|+\sum_{\substack{j\in I\\j>r}}2jL_j+2r(\ell-1)\\
 &\le b+\sum_{j=r+1}^h2jL_j+2r(L_r-1)\\
 &=Q_r-1<Q_r.
\end{aligned}
\]

If $d_{r,\ell}=0$, take $U_{r,\ell}=\emptyset$; then \eqref{eq17} is clearly satisfied. We may therefore assume that $d_{r,\ell}\ne0$.
Since $d_{r,\ell}\in D_{r,Q_r}$, there are pairwise disjoint nonempty sets
$G_1,\ldots,G_{Q_r}\subseteq A$ such that
$$d_{r,\ell}=\sum_{u\in G_i}u,~~1\le |G_i| \le2r ~(1\le i\le Q_r),$$ and
 $G_i\cap G_{i'} =\emptyset~(i\ne i').$
Since $G_1,\ldots,G_{Q_r}$ are pairwise disjoint,
the number of sets meeting $W_{r,\ell}$ is at most $|W_{r,\ell}|$. Indeed,
\[
\begin{aligned}
 \left|\{i:G_i\cap W_{r,\ell}\ne\emptyset\}\right|
 \le\sum_{i=1}^{Q_r}|G_i\cap W_{r,\ell}|
 =\left|\bigcup_{i=1}^{Q_r}(G_i\cap W_{r,\ell})\right|
 \le|W_{r,\ell}|\le Q_r-1.
\end{aligned}
\]
Hence,
\[
\begin{aligned}
 \left|\{i:G_i\cap W_{r,\ell}=\emptyset\}\right|
 =Q_r-\left|\{i:G_i\cap W_{r,\ell}\ne\emptyset\}\right|
 \ge Q_r-(Q_r-1)=1.
\end{aligned}
\]
Thus there is an index $i_0\in\{1,2,\ldots,Q_r\}$ such that $G_{i_0}\cap W_{r,\ell}=\emptyset$. Set $U_{r,\ell}=G_{i_0}$. Then
\[d_{r,\ell}=\sum_{u\in U_{r,\ell}}u,~~|U_{r,\ell}|\le2r,~~U_{r,\ell}\cap W_{r,\ell}=\emptyset.\]
Moreover, by \eqref{eq18}, $F\subseteq W_{r,\ell}$ and $U_{j,m}\subseteq W_{r,\ell}$
whenever $(j,m)\prec(r,\ell)$. Since $U_{r,\ell}\cap W_{r,\ell}=\emptyset$, we have
$U_{r,\ell}\cap F=\emptyset$ and
$U_{r,\ell}\cap U_{j,m}=\emptyset$
whenever $(j,m)\prec(r,\ell)$.

This proves the case $k=1$ and the induction step for $k>1$.
There are $\sum_{r\in I}L_r$ index pairs in $\mathcal J$. Thus, by induction, we can choose all the required sets so that
\begin{equation}\label{eq19}
 d_{r,\ell} =\sum_{u\in U_{r,\ell}}u,\quad |U_{r,\ell}|\le2r,
 U_{r,\ell}\cap F =\emptyset,\quad
 U_{r,\ell}\cap U_{r',\ell'} =\emptyset~((r,\ell)\ne(r',\ell')).
\end{equation}

Finally, put
\[
 U=\bigcup_{r\in I}\bigcup_{\ell=1}^{L_r}U_{r,\ell}.
\]
By \eqref{eq19}, we have $U\subseteq A\setminus F$, and
\[s=\sum_{r\in I}\sum_{\ell=1}^{L_r}d_{r,\ell}=\sum_{r\in I}\sum_{\ell=1}^{L_r}\sum_{u\in U_{r,\ell}}u=\sum_{u\in U}u,\]
\[|U|=\sum_{r\in I}\sum_{\ell=1}^{L_r}|U_{r,\ell}|\le\sum_{r\in I}2rL_r\le\sum_{r=1}^h2rL_r=K.\]
Thus $s\in R_K(A\setminus F)$. Since $s\in S$, $F\subseteq A$, and
$|F|\le b$ were arbitrary, \eqref{eq16} gives
$$[N_0,\infty)\cap(a+g\Z)
 \subseteq S\subseteq R_K(A\setminus F)$$
for $|F|\le b$.

This completes the proof of Lemma~\ref{lem4}.
\end{proof}

\section{Locally complete implies uniformly locally complete}

Let $h=\ord(A)$. For $q,n\ge1$, put
$$P_q=\{r\in\Z/q\Z:\left|\{a\in A:a\equiv r~(\text{mod}~q)\}\right|=\infty\},$$
$C_q=P_q\cup\{0\}$, and $M_n=n!$.

We first compare the sets for different moduli. If $q\mid Q$, then for each $r\in P_Q$,
\[\left|\{a\in A:a\equiv r\pmod q\}\right|=\infty,\]
and hence
\[r\bmod q\in P_q.\]
Therefore,
\begin{equation}\label{eq23}
 C_Q\bmod q\subseteq C_q\quad(q\mid Q).
\end{equation}
We recall that $jC_q$ denotes the set of sums of $j$ elements of $C_q$, where repetitions
are allowed.

We claim that there are $1\le j\le h$, $u\in S_{h-j}(A)$, $m\ge1$, and
$0\le t<M_m$ such that, for every $n\ge m$,
\begin{equation}\label{eq24}
 \{r\in\Z/M_n\Z:r\equiv t~(\text{mod}~{M_m})\}
 \subseteq u+jC_{M_n}.
\end{equation}
Suppose, for a contradiction, that \eqref{eq24} does not hold. Fix
\[
 1\le j\le h,\quad u\in S_{h-j}(A),\quad
 m\ge1,\quad 0\le t<M_m.
\]
Then, for each fixed choice of $j,u,m,t$,
there are an $n\ge m$ and a residue class $r\bmod M_n$ such that
\[
 r\equiv t~(\text{mod}~{M_m}),\quad r\notin u+jC_{M_n}.
\]
We may also require $n>m$ so that the moduli increase strictly.
Indeed, if the chosen value is $n=m$, then
$t\notin u+jC_{M_m}$ and
$$(u+jC_{M_{m+1}})\bmod M_m\subseteq u+jC_{M_m}.$$
This inclusion follows from \eqref{eq23}. Hence,
\[
 \{r'\bmod M_{m+1}:r'\equiv t~(\text{mod}~{M_m})\}
 \cap(u+jC_{M_{m+1}})=\emptyset.
\]
In this case, we may instead take $n=m+1$ and any residue class
$r'\equiv t\pmod{M_m}$.

List all pairs
with $1\le j\le h$ and $u\in S_{h-j}(A)$ in a sequence, repeating pairs if necessary,
\[
 (j_1,u_1),(j_2,u_2),\ldots
\]
Set
\[
 n_0=1,~r_0=0~(\text{mod}~{M_{n_0}}).
\]
At the $k$th step, take $m,t,j,u$ in the preceding paragraph to be
$n_{k-1},r_{k-1},j_k,u_k$, respectively. Since the assumption holds for every choice of parameters,
we can choose, one after another,
\begin{equation}\label{eq25}
 \begin{aligned}
 n_k >n_{k-1}, ~
 r_k \in\Z/M_{n_k}\Z, ~
 r_k \equiv r_{k-1}~(\text{mod}~{M_{n_{k-1}}}), ~
 r_k \notin u_k+j_kC_{M_{n_k}}~(k\ge1).
 \end{aligned}
\end{equation}

Since $M_n\mid M_{n+1}$, \eqref{eq25}
implies that
\[
 r_\ell\equiv r_k~(\text{mod}~{M_{n_k}})\quad(\ell\ge k).
\]
Thus, for any $n\ge1$, choose $k$ such that $n_k\ge n$, and define
\[
 x_n=r_k\bmod M_n.
\]
If $\ell\ge k$ and $n_k\ge n$,
then, since $M_n\mid M_{n_k}$, the preceding congruence also gives
\[
 r_\ell\equiv r_k\pmod{M_n}.
\]
Thus $x_n$ is independent of the choice of $k$. Moreover, since $M_n\mid M_{n+1}$,
\begin{equation}\label{eq26}
 x_{n+1}\equiv x_n~(\text{mod}~{M_n}),\quad
 x_{n_k}=r_k\quad(k\ge1).
\end{equation}

Since $h=\ord(A)<\infty$, there is an $N_A\ge1$ such that
\[
 [N_A,\infty)\cap\Z\subseteq S_h(A).
\]
Every residue class modulo $M_n$ contains arbitrarily large positive integers.
Therefore, we can choose $s_n\ge\max\{n,N_A\}$ such that
$s_n\equiv x_n\pmod{M_n}$.
Since $A$ is an asymptotic basis, $s_n\in S_h(A)$.
Choose a representation with at most $h$ summands and add zeros if necessary. We obtain
$s_n\equiv x_n\pmod{M_n}$ for $s_n\ge n$, and
$$s_n=a_{1,n}+\cdots+a_{h,n},\quad
 a_{i,n}\in A\cup\{0\}.$$

We now select subsequences for only finitely many coordinates.
For $i=1,\ldots,h$, consider $(a_{i,n})$ successively. If some value
$b_i$ occurs infinitely often, keep those indices so that this coordinate is constantly equal to $b_i$.
If every value occurs only finitely often, then for any $T\ge0$,
only finitely many remaining indices satisfy $a_{i,n}\le T$.
Thus this coordinate tends to infinity.
After these $h$ selections, we obtain a strictly increasing sequence $m_1<m_2<\cdots$ and
a set $J\subseteq\{1,\ldots,h\}$ such that
\begin{equation}\label{eq27}
a_{i,m_k}=b_i~(i\notin J,\ k\ge1),\quad
 a_{i,m_k}\rightarrow\infty~(i\in J).
\end{equation}
We must have $J\ne\emptyset$. Otherwise all coordinates are fixed, and
for every $k\ge1$,
\[
 m_k\le s_{m_k}=\sum_{i=1}^h b_i,
\]
which contradicts $m_k\to\infty$. Put
\[
 j=|J|,\quad u=\sum_{i\notin J}b_i,
 \quad 1\le j\le h,\quad u\in S_{h-j}(A).
\]

We next fix a single modulus $M_v$ and prove that
$x_v\in u+jC_{M_v}$. Fix any $v\ge1$ and consider the vector
\[
 (a_{i,m_k}\bmod M_v)_{i\in J}
\]
which has at most $M_v^j$ possible values. By the infinite pigeonhole principle, there are an infinite set of indices
$I_v\subseteq\{k:m_k\ge v\}$ and residue classes $c_{i,v}\bmod M_v$ such that
\[
 |I_v|=\infty,~~ a_{i,m_k}\equiv c_{i,v}~(\text{mod}~{M_v})~~(i\in J,\ k\in I_v).
\]
Since $a_{i,m_k}\to\infty$ for $i\in J$,
the values $a_{i,m_k}$ along this infinite subsequence must include infinitely many distinct elements.
Thus, for each $i\in J$, \eqref{eq27} gives
\[
 \begin{aligned}
 \left|\{a_{i,m_k}:k\in I_v\}\right| =\infty,~~
 \left|\{a\in A:a\equiv c_{i,v}~(\text{mod}~{M_v})\}\right| =\infty.
 \end{aligned}
\]
It follows that
\[c_{i,v}\in P_{M_v}\subseteq C_{M_v}.\]
On the other hand, for every $k\in I_v$, since $m_k\ge v$,
\eqref{eq26}
and $s_{m_k}\equiv x_{m_k}\pmod{M_{m_k}}$ give
\begin{equation}\label{eq29}
 x_v \equiv x_{m_k}\equiv s_{m_k} \equiv\sum_{i\notin J}b_i+\sum_{i\in J}c_{i,v} \equiv u+\sum_{i\in J}c_{i,v}\pmod{M_v},
\end{equation}
Hence, $x_v\in u+jC_{M_v}$.
Since $v$ is arbitrary,
$x_v\in u+jC_{M_v}$ holds for every $v\ge1$.
Choose $k$ such that $(j_k,u_k)=(j,u)$.
By \eqref{eq26} and \eqref{eq29},
\[
 r_k=x_{n_k}\in u+jC_{M_{n_k}}
 =u_k+j_kC_{M_{n_k}},
\]
which contradicts \eqref{eq25}.
Thus \eqref{eq24} holds.

Fix $j,u,m,t$ as in \eqref{eq24}, and put $g=M_m$.
For $g\mid Q$ and $v\in\Z$, write
\[
 T_Q(v)=\{r\in\Z/Q\Z:r\equiv v~(\text{mod}~ g)\}.
\]
Since $Q\mid M_n$ (it suffices to take $n\ge\max\{m,Q\}$),
\eqref{eq23} gives
$$T_Q(t) =T_{M_n}(t)\bmod Q \subseteq(u+jC_{M_n})\bmod Q \subseteq u+jC_Q,$$
and $T_Q(t-u) \subseteq jC_Q.$

Take any $r\in T_Q(0)$. Since $r\equiv0\pmod g$, we have
\[
 r-(g-1)(t-u)\equiv t-u~(\text{mod}~g).
\]
Thus $r$ is the sum of $g-1$ copies of $t-u$ and one copy of
$r-(g-1)(t-u)$, and all these $g$ terms belong to $T_Q(t-u)$.
Conversely, the sum of any $g$ elements of $T_Q(t-u)$ is congruent to zero modulo $g$.
Therefore,
\begin{equation}\label{eq31}
T_Q(0) =\underbrace{T_Q(t-u)+\cdots+T_Q(t-u)}_{g\text{ summands}} \subseteq gjC_Q~~(g\mid Q).
\end{equation}

Since $\Sigma(A)\bmod g=\Z/g\Z$, for each
$v\in\Z/g\Z$, choose a finite set $E_v\subseteq A$ with the smallest possible number of elements such that
\[
 e_v:=\sum_{a\in E_v}a\equiv v\pmod g.
\]
We now show that $|E_v|\le g-1$.
Suppose that $|E_v|\ge g$. Choose distinct $a_1,\ldots,a_g\in E_v$, put $z_0=0$ and
$$z_i=\sum_{k=1}^{i}a_k\quad(1\le i\le g).$$
Since there are $g+1$ integers $z_0,\ldots,z_g$ but only $g$ residue classes modulo $g$,
the pigeonhole principle gives $0\le p<r\le g$ such that $z_p \equiv z_r\pmod g$,
$\sum_{k=p+1}^{r}a_k \equiv0\pmod g$,
$$\sum_{a\in E_v\setminus\{a_{p+1},\ldots,a_r\}}a
 \equiv e_v\equiv v\pmod g,$$
and
$$|E_v\setminus\{a_{p+1},\ldots,a_r\}|
 =|E_v|-(r-p)<|E_v|,$$
contradicting the minimality of $E_v$. Hence
\[
 |E_v|\le g-1\quad(v\in\Z/g\Z).
\]

Fix $q\ge1$ and $w\in\Z/q\Z$.
Choose an integer $x\equiv w\pmod q$ and put
$$v=x\bmod g,\quad Q=\operatorname{lcm}(g,q).$$
Since $e_v\equiv v\equiv x\pmod g$,
we have $x-e_v\equiv0\pmod g$.
Since $g\mid Q$, \eqref{eq31} gives
\[
 (x-e_v)\bmod Q\in T_Q(0)\subseteq gjC_Q.
\]
Thus there are $c_1,\ldots,c_{gj}\in C_Q$ such that
\[
 x-e_v\equiv\sum_{i=1}^{gj}c_i\pmod Q.
\]
Let
$I=\{1\le i\le gj:c_i\ne0\}$.
Since every nonzero $c_i$ belongs to $P_Q$, for each $i\in I$,
\[
 \left|\{a\in A:a\equiv c_i\pmod Q\}\right|=\infty.
\]
Hence we can choose $a_i$ such that
$$a_i \in A\setminus\bigl(E_v\cup\{a_k:k\in I,\ k<i\}\bigr),$$
and
$$a_i\equiv c_i\pmod Q\quad(i\in I).$$
Put $F=E_v\cup\{a_i:i\in I\}$. Then $F\subseteq A$,
$$\sum_{a\in F}a
=e_v+\sum_{i\in I}a_i \equiv e_v+\sum_{i\in I}c_i \equiv e_v+\sum_{i=1}^{gj}c_i
 \equiv x\equiv w\pmod q,$$
$$|F|=|E_v|+|I|\le g-1+gj.$$

Taking $H=g-1+gj$, we obtain
\[
 R_H(A)\bmod q=\Z/q\Z\qquad(q\ge1).
\]
Note that $H$ is independent of $q$.
This completes the proof of Theorem~\ref{thm3}.

\section{Proof of Theorem~\ref{thm2}}

\begin{proof}[Proof of Theorem~\ref{thm2}]The necessity is almost obvious. We prove it briefly. Suppose that $\ordR(A)<\infty$, and
put $j=\ordR(A)$. Then
\[
R_j(A)\sim\mathbb N,\quad R_j(A)\subseteq\Sigma(A).
\]
Hence, for every integer $q\geq1$,
\[
\mathbb Z/q\mathbb Z
=R_j(A)\bmod q
\subseteq\Sigma(A)\bmod q
\subseteq\mathbb Z/q\mathbb Z.
\]
Thus $\Sigma(A)\bmod q=\mathbb Z/q\mathbb Z$ for every $q\geq1$.

Now we start to prove the sufficiency. On applying Theorem \ref{thm3},  there is an integer $H\ge1$ such that
\[
 R_H(A)\bmod q=\Z/q\Z\quad(q\ge1).
\]
Since $A$ has finite ordinary order, choose $h\ge1$ such that $S_h(A)\sim\mathbb N$.
By Lemma~\ref{lem2},
\[
 \ld(R_h(A))>0.
\]
Applying Lemma~\ref{lem4} with $b=H$, we obtain integers $K,g,N_0\ge1$ and $a$
with $0\le a<g$ such that, for every finite set $F\subseteq A$ with $|F|\le H$,
\begin{equation}\label{eq21}
 [N_0,\infty)\cap(a+g\Z)\subseteq R_K(A\setminus F).
\end{equation}
For each $v\in\Z/g\Z$, choose a finite set $E_v\subseteq A$ such that $|E_v|\le H$ and
\[e_v:=\sum_{u\in E_v}u\equiv v-a\pmod g.
\]
Put
\[
 M=N_0+\max_{v\in\Z/g\Z}e_v.
\]
For any $n\ge M$, let $v=n\bmod g$. Then $n-e_v \ge M-e_v\ge N_0$ and
$$n-e_v\equiv v-(v-a)\equiv a\pmod g.$$
Taking $F=E_v$ in \eqref{eq21} (where $|E_v|\le H$),
we obtain a finite set $U_n\subseteq A\setminus E_v$ such that
\[|U_n|\le K,\quad \sum_{u\in U_n}u=n-e_v,\quad U_n\cap E_v=\emptyset.\]
Hence,
\[ n=\sum_{u\in U_n\cup E_v}u,\quad
 |U_n\cup E_v|\le K+H.\]

Since $E_v\cap U_n=\emptyset$, the summands in this representation are pairwise distinct.
Since $n\ge M$ is arbitrary, we finally obtain $[M,\infty)\cap\N\subseteq R_{K+H}(A)$, and
$$\ordR(A)\le K+H<\infty.$$

This completes the proof of Theorem~\ref{thm2}. \end{proof}

\section{Proof of Theorem~\ref{thm1}}

\begin{proof}[Proof of Theorem~\ref{thm1}]For each integer $N\ge1$, put
\[
 A_N=A\cap[N,\infty).
\]
Since $A\setminus A_N=A\cap[0,N-1]$ is finite,
the hypothesis gives an integer $h_N\ge1$ such that $S_{h_N}(A_N)$
contains all sufficiently large integers. We first show that
\begin{equation}\label{eq33}
 \gcd(A_N)=1\quad(N\ge1).
\end{equation}
Indeed, if $d=\gcd(A_N)\ge2$, then $A_N\subseteq d\N$,~$S_k(A_N)\subseteq d\mathbb N~(k\ge1)$,
which contradicts the fact that $A_N$ is an asymptotic basis of finite order. Thus \eqref{eq33} holds.

Now fix an integer $q\ge1$, and put
\[
 \mathcal P_q=\{u\in\{0,1,\ldots,q-1\}:|A\cap(u+q\Z)|=\infty\}.
\]
Since $A$ is infinite,
we may write $\mathcal P_q$ as
\[
 \mathcal P_q=\{u_1,\ldots,u_t\},\quad t\ge1.
\]
Put
\[
 d=\gcd(q,u_1,\ldots,u_t).
\]
We use \eqref{eq33} to show that $d=1$.
Suppose that $d\ge2$. Then $q\equiv0\pmod d$,
$$u_i\equiv0\pmod d\quad(1\le i\le t),$$
and
$$\left|\{a\in A:a\bmod q\notin \mathcal P_q\}\right|<\infty.$$
Since each residue class not belonging to $\mathcal P_q$ contains only finitely many elements of $A$,
we may choose $N$ sufficiently large so that every element of $A_N$ lies in one of the residue classes listed in $\mathcal P_q$.
Since $d\mid q$ and $d\mid u_i$, we obtain $A_N\subseteq d\N$ and
$d\mid\gcd(A_N)=1$,
a contradiction. Therefore,
\[
 \gcd(q,u_1,\ldots,u_t)=1.
\]
Thus there are integers $z_0,z_1,\ldots,z_t$ such that
\[
 1=z_0q+\sum_{i=1}^{t}z_i u_i.
\]
Fix any $v\in\{0,1,\ldots,q-1\}$.
Choose nonnegative integers $k_i$ such that
\[
 k_i\equiv vz_i\pmod q\quad(1\le i\le t).
\]
Then
\[
 v\equiv\sum_{i=1}^{t}k_i u_i\pmod q.
\]
Since every $u_i\in \mathcal P_q$ occurs infinitely often in $A$,
for each $1\le j\le k_i$,
we can choose pairwise distinct elements $a_{i,j}\in A$ such that
\[
 a_{i,j}\equiv u_i\pmod q.
\]
It follows that
\[\sum_{i=1}^{t}\sum_{j=1}^{k_i}a_{i,j}\equiv\sum_{i=1}^t k_i u_i\equiv v\pmod q.\]
Noting that
 \[\sum_{i=1}^{t}\sum_{j=1}^{k_i}a_{i,j}\in\Sigma(A)\]
and $q\ge1$ was arbitrary, we have
\[
 \Sigma(A)\bmod q=\Z/q\Z\quad(q\ge1).
\]

Finally, applying Theorem~\ref{thm2} directly, we obtain
\[
 \ordR(A)<\infty.
\]
This completes the proof of Theorem~\ref{thm1}.\end{proof}

\end{document}